\documentclass[11pt]{article}

\usepackage{graphicx}
\usepackage{url}
\usepackage{color}
\definecolor{maroon}{rgb}{.69,.188,.376}
\definecolor{darkgreen}{rgb}{0,.5,0}
\definecolor{darkblue}{rgb}{0,0,.5}
\definecolor{magenta}{rgb}{1,0,1}

\usepackage[mathscr]{euscript}		

\usepackage{dsfont}

\usepackage{psfrag}			

\usepackage[colorlinks=true]{hyperref}
\hypersetup{pdftex, colorlinks=true, linkcolor=maroon, citecolor=maroon,
  filecolor=blue,urlcolor=blue}

\usepackage{amsmath}
\usepackage{amsthm}

\usepackage{amssymb}

\usepackage{caption} 

\usepackage{anysize}

\usepackage{enumerate}
\usepackage{enumitem}

\usepackage[top=.95in, bottom = 1 in, left=1in, right = .6in]{geometry}

\newcommand{\ind}{\mathds}

\newcommand{\Z}{\ensuremath{\mathbb{Z}}}
\newcommand{\N}{\ensuremath{\mathbb{N}}}
\newcommand{\R}{\ensuremath{\mathbb{R}}}

\newcommand{\E}{\ensuremath{\mathbb{E}}}
\renewcommand{\P}{\ensuremath{\mathbb{P}}}

\newtheorem{theorem}{Theorem}[section]

\newtheorem{remark}[]{Remark}

\newtheorem{definition}[theorem]{Definition}

\numberwithin{equation}{section}

\definecolor{Red}{rgb}{1,0,0}
\definecolor{Blue}{rgb}{0,0,1}
\definecolor{Olive}{rgb}{0.41,0.55,0.13}
\definecolor{Yarok}{rgb}{0,0.5,0}
\definecolor{Green}{rgb}{0,1,0}
\definecolor{MGreen}{rgb}{0,0.8,0}
\definecolor{DGreen}{rgb}{0,0.55,0}
\definecolor{Yellow}{rgb}{1,1,0}
\definecolor{Cyan}{rgb}{0,1,1}
\definecolor{Magenta}{rgb}{1,0,1}
\definecolor{Orange}{rgb}{1,.5,0}
\definecolor{Violet}{rgb}{.5,0,.5}
\definecolor{Purple}{rgb}{.75,0,.25}
\definecolor{Brown}{rgb}{.75,.5,.25}
\definecolor{Grey}{rgb}{.7,.7,.7}
\definecolor{Black}{rgb}{0,0,0}

\newcommand{\ignore}[1]{{}}

\date{\today}
\begin{document}
\baselineskip=14pt

\title{Random Walk Among Mobile/Immobile Traps: A Short Review}

\author{
Siva Athreya%
  \thanks{8th Mile Mysore Road, Indian Statistical Institute,
         Bangalore 560059, India.
    Email: \url{athreya@isibang.ac.in}}
  \and
  Alexander Drewitz%
  \thanks{Universit\"at zu K\"oln,
Mathematisches Institut,
Weyertal 86--90,
50931 K\"oln, Germany.
    Email: \url{drewitz@math.uni-koeln.de}}
  \and
 Rongfeng Sun%
  \thanks{Department of Mathematics, National University of Singapore,
S17, 10 Lower Kent Ridge Road
Singapore, 119076.
    Email: \url{matsr@nus.edu.sg}}
}

\maketitle

\centerline{\em In celebration of Chuck's 70th Birthday}

\bigskip

\begin{abstract}
There have been extensive studies of a random walk among a field of immobile traps (or obstacles), where one is interested in the probability of survival as well as the law of the random walk conditioned on its survival up to time $t$. In contrast, very little is known when the traps are mobile. We will briefly review the literature on the trapping problem with immobile traps, and then review some recent results on a model with mobile traps, where the traps are represented by a Poisson system of independent random walks on $\Z^d$. Some open questions will be given at the end.

\end{abstract}

\noindent {\em AMS 2010 Subject Classification:} 60K37, 60K35, 82C22.\\
\noindent {\em Keywords:} trapping problem, parabolic Anderson model, random walk in random potential.

\section{Introduction}

The trapping problem, where particles diffuse in space with randomly located traps, has been studied extensively in the physics and mathematics literature. We refer the reader to the review article \cite{HW94}, which explains in detail the background for the trapping problem and some early results. Here we focus on a single particle diffusing on $\Z^d$ according to a random walk, with randomly located traps that may or may not be mobile. When the particle meets a trap, it is killed at a fixed rate $\gamma\in (0,\infty]$.

More precisely, let $X:=(X(t))_{t\geq 0}$ be a random walk on $\Z^d$ with jump rate $\kappa \geq 0$.
Let $\xi:= (\xi(t,\cdot))_{t \ge 0}$ be a continuous time Markov process with values in $[0,\infty]^{\Z^d}$, which determines the trapping potential. Let $\gamma\in (0,\infty]$. If $X(t)=x$, then it is killed at rate $\gamma \xi(t,x)$. As we will further detail below, there are two fundamentally different regimes: The setting of {\em immobile} or {\em static traps}, where $\xi(t,\cdot)$ is constant in $t$ (corresponding to the traps being realized at time $0$ and not evolving in time), and the setting of {\em mobile traps}, where $\xi(t,\cdot)$ depends non-trivially on $t$. One can also consider the continuum model where $X$ is a Brownian motion on $\R^d$ and $\xi(t, \cdot)\in [0,\infty]^{\R^d}$, which we will do when reviewing some classic results.

Denote by $\P^\xi$ and $\E^\xi$  the probability of and expectation with respect to $\xi,$
and similarly  by $\P^X_x$ and $\E^X_x$  the probability of and expectation with respect to $X$ when starting at $x \in \Z^d.$
The above model gives rise to the following quantities of interest.

\begin{definition}[Survival Probabilities]\label{D:surprob}
Conditional on the realization of $\xi,$ the {\em quenched survival probability of $X$ up to time $t$} is defined by
\begin{equation}\label{eq:qsurprob}
Z^{\xi}_{\gamma, t} := \E^X_0\Big[\exp\Big\{-\gamma \int_0^t \xi(s,X(s)) \, {\rm d}s\Big\}\Big].
\end{equation}
The {\em annealed survival probability up to time $t$} is defined by averaging over the trap configuration:
\begin{equation}\label{eq:asurprob}
Z_{\gamma, t} := \E^\xi \big[ Z^{\xi}_{\gamma, t} \big] =\E^\xi\Big[ \E^X_0\Big[\exp\Big\{-\gamma \int_0^t \xi(s,X(s)) \, {\rm d}s\Big\}\Big]\Big]
.
\end{equation}
\end{definition}

The above expressions immediately lead to the following definition of path measures, which are the laws of the random walk $X$ conditioned on survival up to time $t$.
\begin{definition}[Path Measures]\label{D:pathmeas}
We call the family of Gibbs measures
\begin{equation} \label{eq:quenchedGibbs}
P_{\gamma, t}^\xi ( X\in \cdot ) := \frac{\E_0^X \Big[  \exp \Big \{ -\gamma \int_0^t \xi(s,X(s)) \, {\rm d}s \Big\} \ind{1}_{X \in \cdot} \Big] }
{   Z^\xi_{\gamma, t}}, \quad t \ge 0,
\end{equation}
on the space of c\`adl\`ag paths $D([0,t], \Z^d)$ from $[0,t]$ to $\Z^d$ the {\em quenched path measures}.

Similarly, the family
\begin{equation} \label{eq:annealedGibbs}
P_{\gamma, t} ( X\in \cdot ) := \frac{\E_0^X \Big[ \E^\xi  \Big[ \exp \Big \{ -\gamma \int_0^t \xi(s,X(s)) \, {\rm d}s \Big\} \Big] \ind{1}_{X \in \cdot} \Big] }
{  \E^\xi[Z^\xi_{\gamma, t}]}, \quad t \ge 0,
\end{equation}
 will be called the {\em annealed path measures}.
\end{definition}
\noindent
The quenched and annealed survival probabilities and their respective path measures are the key objects of interest for the trapping problem.

The trapping problem is closely linked to the so-called parabolic Anderson model (PAM), which is the solution of the following parabolic equation with random potential $\xi$:
\begin{equation} \label{eq:PAM}
\left.
\begin{aligned}
\frac{\partial}{\partial t} u(t, x) &= \kappa \Delta u(t, x) - \gamma\, \xi(t, x)\, u(t, x),  \\
u(0,x) &= 1,
\end{aligned}
\right.
\qquad \qquad x\in \Z^d,\ t\geq 0,
\end{equation}
where $\gamma, \kappa$ and $\xi$ are as before, and the discrete Laplacian on $\Z^d$ is given by $$\Delta f(x) = \frac{1}{2d} \sum_{\Vert y-x\Vert =1} (f(y)-f(x)).$$ By the Feynman-Kac formula, the solution $u$ is given by
\begin{equation} \label{feynkac}
u(t,0) = \E^X_0\left[\exp\left\{-\gamma\int_0^t \xi(t-s, X(s)) \, {\rm d}s \right\} \right],
\end{equation}
which differs from $Z^\xi_{\gamma, t}$ in (\ref{eq:qsurprob}) by a time reversal. Note that if $\xi(t,\cdot)$ does not depend on time, then $u(t,0)=Z^\xi_{\gamma, t}$, and more generally, if the law of $(\xi(s,\cdot))_{0\leq s\leq t}$ is invariant under time reversal, then
\begin{equation}
   \label{annealedequiv}
 \E^\xi[u(t,0)] = \E^\xi\E^X_0\!\!\left[\exp\!\left\{\!-\gamma\!\int_0^t\! \xi(t-s, X(s)) \, {\rm d}s\!\right\}\!\right]\! =\! \E^X_0\E^\xi\!\!\left[\exp\!\left\{\!-\gamma\!\int_0^t \!\xi(s, X(s))\, {\rm d}s\!\right\}\! \right] = \E^\xi[Z^\gamma_{t,\xi}].
\end{equation}
Thus the study of the trapping problem is intimately linked to the study of the PAM, especially the case of immobile traps. A comprehensive account can be found in the recent monograph by K\"onig \cite{K16}.

Historically, most studies of the trapping problem have focused on the case of immobile traps, for which we now have a very good understanding, see e.g.\ Sznitman \cite{S98} and  \cite{K16}. Our goal is to review some of the recent results in this direction, where $\xi$ is the occupation field of a Poisson system of independent random walks. These results provide first steps in the investigation, while much remains to be understood.

The trapping problem has connections to many other models of physical and mathematical interest, such as chemical reaction networks, random Schr\"odinger operators and Anderson localization, directed polymers in random environment, self-interacting random walks, branching random walks in random environment, etc. The literature is too vast to be surveyed here. The interested reader can consult \cite{HW94} for motivations from the chemistry and physics literature and some early mathematical results, \cite{S98} for results on Brownian motion among Poisson obstacles, a continuum model with immobile traps, \cite{DGRS12} for an overview of literature on mobile traps, and \cite{K16} for a comprehensive survey on the parabolic Anderson model, focusing mostly on immobile traps on the lattice $\Z^d$.

The rest of the paper is organized as follows. In Section \ref{S:immobile}, we briefly review the literature on the trapping problem with immobile traps, which has been the theme of the monographs \cite{S98, K16}. We then review in Section \ref{S:mobile} some recent results on the case of mobile traps \cite{DGRS12, ADS16}. Lastly in Section \ref{S:open}, we discuss some open questions.

\section{Immobile Traps} \label{S:immobile}

In this section, we briefly review what is known for the trapping problem with immobile traps, where the trapping potential $\xi(t,x)=\xi(x)$ does not depend on time and $(\xi(x))_{x\in \Z^d}$ are i.i.d., which has been the subject of the recent monograph \cite{K16}. We focus here on the lattice setting, although historically, the first comprehensive mathematical results were obtained for the problem of {\em Brownian motion among Poisson obstacles} \cite{DV75, S98}, where the random walk is replaced by a Brownian motion in $\R^d$, and the traps are balls whose centers follow a homogeneous Poisson point process on $\R^d$. As we will explain, the basic tools in the analysis of the trapping problem with i.i.d.\ immobile traps are large deviation theory and spectral techniques.

\subsection{Annealed Asymptotics}

We first consider the annealed survival probability $Z_{\gamma, t}$ (equivalently, $\E^\xi[u(t,0)]$ for the PAM). Using the fact that $(\xi(x))_{x\in\Z^d}$ are i.i.d., we can integrate out $\xi$ to write
\begin{equation}\label{annZ1}
Z_{\gamma, t} = \E^\xi\Big[ \E^X_0\Big[\exp\Big\{-\gamma \int_0^t \xi(X(s)) \, {\rm d}s\Big\}\Big]\Big] = \E^X_0\Big[
e^{\sum_{x\in\Z^d} H(\gamma L_t(x))}\Big],
\end{equation}
where $L_t(x):= \int_0^t 1_{\{X_s=x\}} \, {\rm d}s$ is the local time of $X$, and
\begin{equation}\label{Hgen}
H(t) := \ln \E^\xi[e^{-t\xi(0)}].
\end{equation}
In the special case $\gamma=\infty$ and $\xi(x)$ are i.i.d.\ Bernoulli random variables with $\P(\xi(0)=0)=p$, the so-called {\em Bernoulli trap model with hard traps}, the expression simplifies to
\begin{equation}\label{annZ2}
Z_{\infty, t} =\E^X_0\Big[p^{|{\rm Range}_{s \in [0,t]}(X(s))|}\Big] = \E^X_0\big[e^{|{\rm Range}_{s \in [0,t]}(X(s))|\ln p}\big],
\end{equation}
where ${\rm Range}_{s \in [0,t]}(X(s)):=\{X(s)\in \Z^d: s\in [0,t]\}$ is the range of $X$ by time $t$. The asymptotic analysis of $Z_{\infty, t}$, and its continuum analogue, the Wiener sausage, were carried out by Donsker and Varadhan in a series of celebrated works \cite{DV75, DV79} using large deviation techniques. The basic heuristics is that the random walk chooses to stay within a spatial window of scale $1\ll \alpha_t \ll \sqrt{t}$, and within that window, the random walk occupation time measure realizes an optimal profile. Using the large deviation principle for the random walk occupation time measure on spatial scale $\alpha_t$, one can then optimize over the scale $\alpha_t$ and the occupation time profile to derive a variational representation for the asymptotics of $Z_{\infty, t}$. See e.g.\ \cite[Sec.~4.2]{K16} and the references therein. Here we only sketch how to identify the optimal scale $\alpha_t$.

In representation \eqref{annZ2} there are two competing effects: the exponential factor which becomes large when $|{\rm Range}_t(X)|$ is small, and the probabilistic cost of a random walk having a small range. If, as the Faber-Krahn inequality suggests, we assume that the expectation in \eqref{annZ2}
is attained by ${\rm Range}_t(X)$ being roughly a ball of radius $\alpha_t$, then $e^{|{\rm Range}_t(X)|\ln p} \approx e^{-c_1 \alpha_t^d}$, while $\P^X_0(\sup_{0\leq s\leq t}\Vert X_s\Vert\leq \alpha_t )\approx e^{-c_2 t/\alpha_t^2}$ due to Brownian scaling with $c_1$ and $c_2$ positive constants. Hence,
\begin{equation}
Z_{\infty, t}  = \E^X_0\big[e^{|{\rm Range}_t(X)|\ln p}\big] \approx \exp\big\{-\inf\limits_{1\ll \alpha_t\ll \sqrt{t}} (c_1 \alpha_t^d + c_2 t/\alpha_t^2)\big\} = e^{-c_3 t^{\frac{d}{d+2}}},
\end{equation}
where we find that the optimal scale is $\alpha_t\approx t^{\frac{1}{d+2}}$.

Alternatively, we can first identify the optimal profile for the trapping potential $\xi$, which is to create a clearing free of traps in a ball of radius $\alpha_t$ around the origin, and the random walk in such a potential is then forced to stay within this region. The probability of the first event is of the order $e^{-c_1 \alpha_t^d}$, while the probability of the second event is of the order $e^{-c_2 t/\alpha_t^2}$, which leads to the same optimal choice of $\alpha_t\approx t^{\frac{1}{d+2}}$.

For more general $\xi$ and $\gamma$, one can still analyze \eqref{annZ1} via large deviations of the random walk occupation time measure. Suppose that the random walk is confined to a box $[-R \alpha_t, R\alpha_t]^d$ for some $R>0$ and scale $1\ll \alpha_t\ll t^{1/d}$, and let
$\widetilde L_t(y) := \frac{\alpha_t^d}{t} L_t(\lfloor \alpha_t y\rfloor)$, $y\in [-R, R]^d$,
denote the rescaled occupation time measure. Furthermore, assume that the generating function $H$ in \eqref{Hgen} satisfies the assumption
\begin{equation}\label{assH}
{\rm (H):} \qquad \qquad \lim_{t\uparrow \infty} \frac{H(ty) - yH(t)}{\eta(t)} = \widehat H(y) \neq 0 \ \mbox{ for } y\neq 1,
\end{equation}
for some $\widehat H:(0,\infty)\to \R$ and continuous $\eta:(0,\infty)\to (0,\infty)$ with $\lim_{t\to\infty} \eta(t)/t\in [0,\infty]$. The assumption (H) essentially ensures an appropriate regularity in the right tail of the distribution of $-\xi(0)$. We can then rewrite \eqref{annZ1} as
\begin{align}
Z_{\gamma, t} & = \E^X_0\Big[ \exp\Big\{ \eta\Big(\frac{\gamma t}{\alpha_t^d}\Big) \sum_{x\in\Z^d} \frac{H\Big(\frac{\gamma t}{\alpha_t^d}\widetilde L_t\big(\frac{x}{\alpha_t}\big)\Big) - \widetilde L_t\big(\frac{x}{\alpha_t}\big) H\Big(\frac{\gamma t}{\alpha_t^d}\Big)  }{\eta\Big(\frac{\gamma t}{\alpha_t^d}\Big) }\Big\}\Big] e^{\alpha_t^d H\big(\frac{\gamma t}{\alpha_t^d}\big) } \nonumber \\
& \approx  e^{\alpha_t^d H\big(\frac{\gamma t}{\alpha_t^d}\big)}  \E^X_0\Big[ \exp\Big\{\eta\Big(\frac{\gamma t}{\alpha_t^d}\Big)
\sum_{x\in\Z^d} \widehat H\Big(\widetilde L_t\big(\frac{x}{\alpha_t}\big)\Big)\Big\}
 \Big] \nonumber \\
& \approx e^{\alpha_t^d H\big(\frac{\gamma t}{\alpha_t^d}\big)}  \E^X_0\Big[ \exp\Big\{\eta\Big(\frac{\gamma t}{\alpha_t^d}\Big)\alpha_t^d
\int_{\R^d} \widehat H\big(\widetilde L_t (y) \big)\, {\rm d}y\Big\}\Big] \nonumber \\
& = e^{\alpha_t^d H\big(\frac{\gamma t}{\alpha_t^d}\big)}  \E^X_0\Big[ \exp\Big\{\frac{t}{\alpha_t^2} \int_{\R^d} \widehat H\big(\widetilde L_t (y) \big)\,{\rm d}y\Big\}\Big], \label{annZ3}
\end{align}
where the scale $\alpha_t$ is chosen to satisfy
\begin{equation}\label{etaalphat}
\eta\big(\frac{\gamma t}{\alpha_t^d}\big)\alpha_t^d = \frac{t}{\alpha_t^2},
\end{equation}
so that in \eqref{annZ3}, the exponential term is comparable to the large deviation probability of the random walk confined in a spatial window of scale $\alpha_t$. Similar to the Bernoulli hard trap case, the asymptotics of $Z_{\gamma, t}$ can then be identified using the large deviation principle for the random walk occupation time measure on the scale $\alpha_t$.

Alternatively, we can first identify the optimal profile $\psi(\alpha_t\cdot)$ for the trapping potential $\xi$ on the spatial scale $\alpha_t$. The survival probability of a random walk in such a potential then decays like $e^{\lambda_\psi t/\alpha_t^2}$, where $\lambda_\psi$ is the principal eigenvalue of $\Delta-\gamma \psi$. In principle, the asymptotics of $Z_{\gamma, t}$ can then be identified by applying a large deviation principle for the potential $\xi$ on the spatial scale $\alpha_t$, and then optimizing over $\alpha_t$ and $\psi$. In practice, this approach has been difficult to implement, and the large deviation approach outlined above has been the standard route in the study of the PAM~\cite[Sec.~3.2]{K16}.

We remark that the heuristics above applies when the scale $\alpha_t$ chosen as in \eqref{etaalphat} tends to infinity as $t$ tends to infinity. However, there are also interesting cases where $\alpha_t$ remains bounded or even $\alpha_t=1$. Under assumption (H) in \eqref{assH}, there are in fact four classes of potentials~\cite{HKM06}, each determined by the right tail probability of $-\xi(0)$:

\begin{enumerate}
\item[(a)] potentials with tails heavier than those of the double-exponential distribution;

\item[(b)] double-exponentially distributed potentials, i.e., $\P(-\xi(0)>r) = e^{-e^{r/\rho}}$ for some $\rho\in (0,\infty)$;

\item[(c)] so-called {\em almost bounded} potentials;

\item[(d)] bounded potentials.
\end{enumerate}
What distinguishes the four classes are different scales $\alpha_t$, with $\alpha_t=1$ in case (a), also known as the single peak case; $\alpha_t$ stays bounded in case (b), where the potential follows the double exponential distribution; $1\ll \alpha_t\ll t^\epsilon$ for any $\epsilon>0$ in case (c); and $\alpha_t\to\infty$ faster than some power of $t$ in case (d), which includes in particular the Bernoulli trap model, the discrete analogue of Brownian motion among Poisson obstacles. In the annealed setting, the potential $\xi$ realizes an optimal profile on a so-called {\em intermittent island} of spatial scale $\alpha_t$ centered around the origin, and the walk then stays confined in that island. For further details, see \cite[Chapter 3]{K16} and the references therein.

The heuristics sketched above also suggests what the annealed path measure $P_{\gamma, t}$ should look like, namely, the random walk $X$ should fluctuate on the spatial scale $\alpha_t$. If $\alpha_t\to\infty$, then after diffusively rescaling space-time by $(\alpha_t^{-1}, \alpha_t^{-2})$, we expect the random walk to converge to a Brownian motion $h$-transformed by the principal eigenfunction of $\Delta -\gamma \psi$, where $\psi(\alpha_t\cdot)$ is the optimal profile the potential $\xi$ realizes on the intermittent island of scale $\alpha_t$, and with Dirichlet boundary condition off the intermittent island. In practice, however, identifying the path behavior requires obtaining second-order asymptotics for the annealed survival probability $Z_{\gamma, t}$, and complete results have only been obtained in special cases.

Indeed, for the continuum model of a Brownian motion $(X(s))_{s \in [0,\infty)}$ among Poisson obstacles in $d=1,$ Schmock \cite{S90} (hard obstacles) and Sethuraman \cite{S03} (soft obstacles) have shown that under the annealed path measure $P_{\gamma, t}$, the law of $(t^{-\frac13} X(s\cdot t^\frac23))_{s \in [0,\infty)}$ converges weakly to a mixture of the laws of so-called Brownian taboo processes, i.e., a mixture of Brownian motions conditioned to stay inside a randomly centered interval of length $(\pi^2/\nu)^\frac13$, and the distribution of the random center can also be determined explicitly (here $\nu$ is the Poisson intensity of the traps). In particular, this result implies that typical annealed fluctuations for the Brownian motion among Poisson obstacles are of the order $t^\frac13$ in $d=1$. Similar convergence results have been proved by Sznitman in dimension $2$ \cite{S91}, and as observed by Povel in \cite{P99}, can also be established in dimensions $d\geq 3$, where Brownian motion is confined to a ball with radius of the order $t^\frac{1}{d+2}$. For the Bernoulli trap model in dimension $2$, path confinement to a ball with radius of the order $t^\frac{1}{4}$ has been proved by Bolthausen in \cite{B94}.

In the lattice setting, when the potential $(\xi(x))_{x\in\Z^d}$ is i.i.d.\ with double exponential distribution (class (b) above), the intermittent islands are bounded in size, and it is known that the rescaled occupation time measure $\widetilde L_t(x) =L_t(x)/t$, $x\in\Z^d$,
converges to a deterministic measure with a random shift, and the distribution of the random shift can also be identified explicitly \cite[Theorem 7.2]{K16}. The convergence of the annealed path measure should follow by the same proof techniques, although it does not seem to have been formulated explicitly in the literature. For potentials in class (a) above, where the intermittent island is a single site, the same result should hold, and the quenched setting leads to more interesting results (see e.g.\ the survey \cite{M11}).

\subsection{Quenched Asymptotics} \label{sec:QA}
In the quenched setting, the random walk $X$ must seek out regions in space that are favorable for its survival. The heuristics is that, one may first confine the random walk to a box $\Lambda_t$ of spatial scale $t(\ln t)^2$, and within $\Lambda_t$, there are intermittent islands of an optimal spatial scale $\widetilde \alpha_t$, on which $\xi$ is close to $\min_{x\in \Lambda_t} \xi(x)$ and takes on an optimal deterministic profile that favors the random walk's survival. The random walk then seeks out an optimal intermittent island and stays there until time $t$, where the choice of the island depends on the balance between the cost for the random walk to get there in a short time and the probability of surviving on the island until time $t$.

For instance, for Brownian motion among Poisson obstacles (or the
Bernoulli trap model), the intermittent islands should be balls
containing no obstacles, and it is then easy to see that within
$\Lambda_t$, the largest such intermittent islands should have spatial
scale of the order $\widetilde \alpha_t=(\ln t)^{1/d}$, which
suggests that the quenched survival probability $Z^\xi_{\gamma, t}$
should decay like $e^{-c t/(\ln t)^{2/d}}$. For i.i.d.\ potential
$\xi$ on $\Z^d$ satisfying assumption (H) in \eqref{assH}, the scale
$\widetilde \alpha_t$ for the intermittent islands can be determined
similarly by considering the following:
\begin{enumerate}
\item the large deviation probability of $\xi$ achieving a
non-trivial profile on the intermittent island after suitable
centering and scaling;
\item a balance between this large deviation probability and the
probability that the random walk survives on the intermittent island
until time $t$;
\item and the requirement that there should be of order one
such islands in $\Lambda_t$.
\end{enumerate} It turns out that $\widetilde \alpha_t =
\alpha_{\beta(t)}$, where $\beta(t)$ is another scale satisfying
\begin{equation}
\frac{\beta(t)}{\alpha_{\beta(t)}^2} = \ln (t^d).
\end{equation}
See e.g.~\cite[Sec.~5.1]{K16} for further details. Note that if $\alpha_t= t^{\epsilon}$ for some $\epsilon\in [0, 1/2)$, then $\beta(t)= (d\ln t)^{\frac{1}{1-2\epsilon}}$. This implies that the size of the quenched intermittent islands is of the same order as the size of the annealed intermittent islands when the latter is bounded, and is much smaller when the latter grows to infinity. The heuristic picture above leads to a sharp lower bound on $Z^\xi_{\gamma, t}$. The challenge is to obtain the matching upper bound, as well as to identify finer asymptotics to draw conclusions about the quenched path measure $P^\xi_{\gamma, t}$.

The main tools are spectral techniques. The exponential rate of decay of the quenched survival probability $Z^\xi_{\gamma, t}$ is given by the principal eigenvalue of $\Delta-\gamma \xi$ on the box $\Lambda_t$. Different techniques have been developed to bound the principal eigenvalue. For Brownian motion among Poisson obstacles, Sznitman developed the {\em method of enlargement of obstacles} (MEO), which resulted in the monograph \cite{S98}. The basic idea is to enlarge the obstacles to reduce the combinatorial complexity, without significantly altering the principal eigenvalue of $\Delta-\gamma\xi$ in $\Lambda_t$. We refer to the review \cite{K00} and the monograph \cite{S98} for further details. The MEO was adapted to the lattice setting in \cite{A95}.
In the lattice setting, a different approach has been developed in the context of the PAM. One divides $\Lambda_t$ into microboxes of scale $\widetilde \alpha_t$. The principal Dirichlet eigenvalues of $\Delta-\gamma \xi$ on the microboxes are i.i.d.\ random variables, and spectral domain decomposition techniques allow one to control the principal eigenvalue on $\Lambda_t$ in terms of the principal eigenvalues on the microboxes. This reduces the analysis to an extreme value problem, where the random walk optimizes over the choice of the microbox to go to and the cost of getting there (see \cite[Sec.~4.4, Sec.~6.3]{K16}). We note that there is also a simple way to transfer the annealed asymptotics for the survival probability to the quenched asymptotics using the so-called Liftshitz tail, as shown by Fukushima in \cite{F09}.

As the heuristics suggest, under the quenched measure, the random walk should seek out one of the optimal intermittent islands and stay there until time $t$. This has been partially verified. For Brownian motion among Poisson obstacles, Sznitman~\cite{S98} used the MEO to show that the quenched survival probability decays asymptotically as $e^{-(C+o(1)) t/(\ln t)^{2/d}}$. The expected picture is that there are $t^{o(1)}$ many intermittent islands, with size $t^{o(1)}$ and mutual distance $t^{1-o(1)}$ to each other and to the origin, and the time when the Brownian motion first enters one of the intermittent islands is $o(t)$. This has been rigorously verified in dimension $1$ (see e.g.~\cite[Chapter 6]{S98}). But in higher dimensions, there is still a lack of control on the time of entering the intermittent islands. Much more precise control on the size and number of intermittent islands have recently been obtained by Ding and Xu for the Bernoulli hard trap model~\cite{DX17}: there are at most $(\log n)^a$ many intermittent islands, with size at most $(\log n)^b$ for some $a, b>0$, and their distance to the origin is at least $n/(\log n)^a$. For the PAM on $\Z^d$ with i.i.d.\ potential $(\xi(x))_{x\in\Z^d}$ satisfying assumption (H) in \eqref{assH}, sharp asymptotics for the quenched survival probability has been obtained in all cases (see \cite[Theorem 5.1]{K16} and the references therein). When the trapping potential $(\xi(x))_{x\in\Z^d}$ are i.i.d.\ with double exponential distribution, the eigenvalue order statistics has been successfully analyzed in~\cite{BK16}, the quenched path measure has been shown to concentrate on a single island at a distance of the order $t/(\ln t \ln\ln\ln t)$ from the origin, and the scaling limit of the
quenched path measure has been obtained (see \cite{BKS16} or \cite[Theorem~6.5]{K16}). Similar results have been obtained for potentials with
heavier tails, for which the intermittent islands consists of single sites (see \cite[Sec.~6.4]{K16} or \cite{M11}).

\section{Mobile Traps}\label{S:mobile}

There have been few mathematical results on the trapping problem with mobile traps. Redig \cite{R94} considered trapping potential $\xi$ generated by a reversible Markov process, such as a Poisson field of independent random walks or the symmetric exclusion process in equilibrium, and he obtained exponential upper bounds on the annealed survival probability using spectral techniques for the process of traps viewed from the random walk. When $\xi$ is generated by a single mobile trap, i.e., $\xi(t,x)= \delta_x(Y_t)$ for a simple symmetric random walk $Y$ on $\Z^d$, Schnitzler and Wolff~\cite{SW12} have computed the asymptotics of the decay of the annealed survival probability via explicit calculations.
The large deviation and spectral techniques outlined in Section \ref{S:immobile} for immobile traps largely fail for mobile traps, and new techniques need to be developed.

Recently, we investigated further the model where $\xi$ is generated by a Poisson field of independent random walks, previously considered in \cite{R94}. More precisely, given $\nu > 0$, let $(N_y)_{y \in \Z^d}$ be a family of i.i.d.\ Poisson random variables
with mean $\nu$. We then start a family of independent random walks
$(Y^{j,y})_{y \in \Z^d, \; 1 \le j \le N_y}$ on $\Z^d$, each with jump
rate $\rho \ge 0$ and $Y^{j,y}:=(Y^{j,y}_t)_{t\geq 0}$ denotes the
path of the $j$-th trap starting from $y$ at time $0$. For $t \ge 0$
and $x \in \Z^d$, we then define
\begin{equation}\label{xidef}
\xi(t,x) := \sum_{y \in \Z^d, \; 1 \le j \le N_y} \delta_x (Y^{j,y}_t)
\end{equation}
which counts the number of traps at site $x$ at time $t.$ When $\rho=0$, we recover the case of immobile traps.

For trapping potential $\xi$ defined as in \eqref{xidef}, there are no systematic tools except the Poisson structure of the traps. The only known results so far are: (1) The quenched survival probability decays at a well-defined exponential rate in all dimensions, while the annealed survival probability decays sub-exponentially in dimensions $d=1$ and $2$ and exponentially in $d\geq 3$~\cite{DGRS12}; (2) The random walk is sub-diffusive under the annealed path measure in dimension $d=1$~\cite{ADS16}. We will review these results in Sections \ref{S:survprob} and \ref{S:pathmeas} below. We will assume for the rest of this section that $\xi$ is defined as in \eqref{xidef}.

We remark that in the physics literature, recent studies of the trapping problem with mobile traps have focused on the annealed survival probability~\cite{MOBC03, MOBC04}, and the particle and trap motion may also be sub-diffusive (see e.g.~\cite{YOLBK08, BAY09}, and also \cite{CS14} for some mathematical results). In the recent mathematics literature, continuum analogues of the results in \cite{DGRS12} on the survival probability were also obtained in \cite{PSSS13}, where the traps move as a Poisson collection of independent Brownian motions in $\R^d$. Further extensions to L\'evy trap motion were carried out in \cite{DSS14}.

\subsection{Decay of survival probabilities}\label{S:survprob}
The precise rate of decay of the annealed and quenched survival probabilities were determined in \cite{DGRS12}. Recall that $\gamma$, $\kappa$, $\rho$ and $\nu$ are respectively the rate of killing per trap, the jump rate of the random walk $X$, the jump rate of the traps, and the density of the traps.

\begin{theorem}\label{T:quenched}{\bf [Quenched survival probability]}
Assume that $d\geq 1$, $\gamma>0$, $\kappa\geq 0$, $\rho>0$ and $\nu>0$, and the traps as well as the walk $X$ follow independent simple symmetric random walks. Then there exists $\lambda^{\rm q}_{d,\gamma,\kappa,\rho, \nu}$ deterministic  such that $\P^\xi$-a.s.,
\begin{equation}\label{quenched}
Z^\xi_{\gamma, t} = \exp\big\{-\lambda^{\rm q}_{d,\gamma,\kappa,\rho, \nu}\,t(1+o(1))\big\} \quad  \mbox{as } t\to\infty.
\end{equation}
Furthermore, $0<\lambda^{\rm q}_{d,\gamma,\kappa,\rho, \nu}\leq \gamma\nu+\kappa$.
\end{theorem}
\noindent
The proof of Theorem \ref{T:quenched} is based on the sub-additive ergodic theorem~\cite{DGRS12}.

\begin{theorem}\label{T:annealed}{\bf [Annealed survival probability]}
Assume that $\gamma\in (0,\infty]$, $\kappa\geq 0$, $\rho>0$ and $\nu>0$, and the traps as well as the walk $X$ follow independent simple symmetric random walks. Then
\begin{equation}\label{annealed}
Z_{\gamma, t} = \E^\xi[Z^\xi_{\gamma, t}]\ =\ \left\{
\begin{aligned}
\exp\Big\{-\nu \sqrt{\frac{8\rho t}{\pi}}(1+o(1))\Big\}, & & \qquad d=1, \\
\exp\Big\{-\nu\pi\rho \frac{t}{\ln t}(1+o(1))\Big\}, & & \qquad d=2,  \\
\exp\Big\{-\lambda^{\rm a}_{d,\gamma,\kappa,\rho, \nu}\, t(1+o(1))\Big\}, & & \qquad d\geq 3.
\end{aligned}
\right.
\end{equation}
Furthermore, 
$$
\lambda^{\rm a}_{d,\gamma,\kappa,\rho, \nu}\geq \lambda^{\rm a}_{d,\gamma,0,\rho, \nu}=\frac{\nu\gamma}{1+\frac{\gamma G_d(0)}{\rho}},
$$ 
where $G_d(0):=\int_0^\infty p_t(0)\, {\rm d}t$ is the Green function of a simple symmetric random walk on $\Z^d$ with jump rate $1$ and transition kernel $p_t(\cdot)$.
\end{theorem}
\noindent
Note that in dimensions $d=1$ and $2$, the annealed survival probability decays sub-exponentially, and the pre-factor in front of the decay rate is independent of $\gamma\in (0,\infty]$ and $\kappa\geq 0$.

  Although Theorems~\ref{T:quenched} and \ref{T:annealed} were proved in \cite{DGRS12} for traps and $X$ following simple symmetric random walks, they can be easily extended to more general symmetric random walks with mean zero and finite variance.

We sketch below the proof of Theorem \ref{T:annealed} and the main tools used in \cite{DGRS12}. The first step is to integrate out the Poisson random field $\xi$ and derive suitable representations for the annealed survival probability $Z_{\gamma, t}$.
\medskip

\noindent
{\bf Representation of $Z_{\gamma, t}$ in terms of random walk range:} Given the Poisson field $\xi$ defined as in \eqref{xidef}, we can integrate out $\xi$ to obtain
\begin{equation} \label{urep}
Z_{\gamma, t}= \E^\xi[Z^\xi_{\gamma, t}] =\E^X_0\E^\xi\!\left[\exp\!\left\{\!-\gamma\!\int_0^t \!\!\xi(s, X(s))\,{\rm d}s\!\right\}\! \right]\! =\E^X_0\Big[ \exp \Big\{\nu\!\! \sum_{y\in\Z^d} (v_X(t,y)-1)\Big\}\Big],
\end{equation}
where conditional on $X$,
\begin{equation} \label{vrep1}
v_X(t, y) = \E^Y_y\left[\exp \left\{ - \gamma\int_0^t \delta_0(Y(s)-X(s))\,{\rm d}s\right\}   \right]
\end{equation}
with $\E^Y_y[\cdot]$ denoting expectation with respect to the motion of a trap starting at $y$.

When $\gamma=\infty$, we can interpret
$$
v_X(t, y) = \P^Y_y\big(Y(s)\neq X(s) \, \forall\, s\in [0,t]\big),
$$
which leads to
\begin{equation} \label{urep'}
\begin{aligned}
Z_{\infty, t} & =\E^X_0\Big[ \exp \Big\{- \nu\!\! \sum_{y\in\Z^d} \P^Y_y\big(Y(s)-X(s)=0 \mbox{ for some } s\in [0,t]\big) \Big\}\Big] \\
& = \E^X_0\Big[ \exp \Big\{- \nu\!\! \sum_{y\in\Z^d} \P^Y_0\big(Y(s)-X(s)=-y \mbox{ for some } s\in [0,t]\big) \Big\}\Big] \\
& = \E^X_0\Big[ \exp \Big\{- \nu \E^Y_0\big[\big|{\rm Range}_{s\in [0,t]}(Y(s)-X(s))\big|\big]\Big\}\Big],
\end{aligned}
\end{equation}
where we recall
\begin{equation}\label{range1}
{\rm Range}_{s\in [0,t]}(Y(s)-X(s)) = \{Y(s)-X(s) \in \Z \, :\, s\in [0, t]\}.
\end{equation}
from \eqref{annZ2}.

When $\gamma \in (0,\infty)$, we can give a similar representation of $Z_{\gamma, t}$ in terms of the range of $Y-X$. More precisely, let
${\cal T}:=\{T_1, T_2, \ldots\}\subset [0,\infty)$ be an independent Poisson point process on $[0,\infty)$ with intensity $\gamma$, and define
\begin{equation}\label{range2}
{\rm SoftRange}_{s\in [0,t]}(Y(s)-X(s)) := \{Y(T_k)-X(T_k): k\in \N, T_k\in [0,t]\}.
\end{equation}
Then we have
\begin{equation}\label{Zrep2}
Z_{\gamma, t} =  \E^X_0\Big[ \exp \Big\{- \nu \E^{Y, \cal T}_0\big[\big|{\rm SoftRange}_{s\in [0,t]}(Y(s)-X(s))\big|\big]\Big\}\Big],
\end{equation}
where $\E^{Y, \cal T}_0$ denotes expectation with respect to both $Y$ and the Poisson point process $\cal T$.
\bigskip

\noindent
{\bf Alternative representation of $Z_{\gamma, t}$:}
We now derive an alternative representation of $Z_{\gamma, t}$ which will turn out useful below e.g.\@ for showing the existence of an asymptotic exponential decay of the survival probability.
For this purpose, note that the time-reversed process $\widetilde \xi(s, \cdot):= \xi(t-s,\cdot)$, $s\in [0,t]$, is also the occupation field of a Poisson system of random walks, where each walk follows the law of $\widetilde Y:=-Y$, which is the same as that of $Y$ by symmetry. We can then write
\begin{equation}
Z_{\gamma, t}=\E^X_0\E^\xi\!\left[\exp\!\left\{\!-\gamma\!\int_0^t \!\!\widetilde\xi(t-s, X(s))\,{\rm d}s\!\right\}\! \right]\! =\E^X_0\Big[ \exp \Big\{\nu\!\! \sum_{y\in\Z^d} (\widetilde v_X(t,y)-1)\Big\}\Big],
\end{equation}
where
\begin{equation}
\widetilde v_X(t, y) = \E^{\widetilde Y}_y\left[\exp \left\{ - \gamma\int_0^t \delta_0(\widetilde Y(t-s)-X(s))\,{\rm d}s\right\}   \right].
\end{equation}
Let $\widetilde L$ denote the generator of $\widetilde Y$. Then by the Feynman-Kac formula, $(\widetilde v_X(t, y))_{t\geq 0,y\in\Z^d}$  solves the equation
\begin{equation} \label{vX}
\begin{aligned}
\frac{\partial}{\partial t} \widetilde v_X(t, y) & = \widetilde L \widetilde v_X(t, y) - \gamma \delta_{X(t)}(y)\, \widetilde v_X(t, y), \\
\widetilde v_X(0,\cdot) & \equiv 1,
\end{aligned}
\qquad y\in\Z^d, \ t\geq 0,
\end{equation}
which implies that $\Sigma_X(t) := \sum_{y \in \Z^d}(\widetilde v_X(t, y) -1)$ is the solution of the equation
\begin{equation}
\begin{aligned}
\frac{\rm d}{{\rm d}t}\Sigma_X(t)  & = - \gamma \widetilde v_X(t, X(t)), \\
\Sigma_X(0) & = 0.
\end{aligned}
\end{equation}
Hence, $\Sigma_X(t) = -\gamma \int_0^t \widetilde v_X(s, X(s))\, {\rm d}s$, which leads to the alternative representation
\begin{equation} \label{urep2}
Z_{\gamma, t} = \E^X_0\left[ \exp \left\{-\nu\gamma \int_0^t \widetilde v_X(s, X(s))\, {\rm d}s \right\}\right].
\end{equation}
\bigskip

\noindent
{\bf Existence of $\lim_{t\to\infty} \frac{1}{t}\ln Z_{\gamma, t}$:} From the representation \eqref{urep2}, it is easily seen that
$\ln Z_{\gamma, t}$ is super-additive, which implies the existence of $\lim_{t\to\infty} \frac{1}{t}\ln Z_{\gamma, t}$. Indeed, for
$t_1, t_2>0$,
\begin{eqnarray}
Z_{\gamma, t_1+t_2} &=& \E^X_0\left[ \exp \left\{-\nu\gamma \int_0^{t_1} \widetilde v_X(s, X(s))\, {\rm d}s \right\}\exp \left\{-\nu\gamma \int_{t_1}^{t_1+t_2} \widetilde v_X(s, X(s))\, {\rm d}s \right\}\right] \nonumber\\
&\geq& \E^X_0\left[ \exp \left\{-\nu\gamma \int_0^{t_1} \widetilde v_X(s, X(s))\, {\rm d}s \right\}\exp \left\{-\nu\gamma \int_{0}^{t_2} \widetilde v_{\theta_{t_1}X}(s, (\theta_{t_1}X)(s))\, {\rm d}s \right\}\right] \nonumber\\
&=& Z_{\gamma, t_1} Z_{\gamma, t_2}, \label{subadd}
\end{eqnarray}
where $\theta_{t_1}X:=((\theta_{t_1}X)(s))_{s\geq 0}=(X(t_1+s)-X(t_1))_{s\geq 0}$, we used the independence of $(X(s))_{0\leq s\leq t_1}$ and $((\theta_{t_1}X)(s))_{0\leq s\leq t_2}$, and the fact that for $s>t_1$,
\begin{eqnarray*}
\widetilde v_X(s, X(s)) &=& \E^{\widetilde Y}_{X(s)}\left[\exp \left\{ - \gamma\int_0^s \delta_0(\widetilde Y(r)-X(s-r))\, {\rm d}r\right\}   \right] \\
&\leq& \E^{\widetilde Y}_{X(s)}\left[\exp \left\{ - \gamma\int_0^{s-t_1} \delta_0(\widetilde Y(r)-X(s-r))\, {\rm d}r\right\}   \right] = \widetilde v_{\theta_{t_1}X}(s-t_1, (\theta_{t_1}X)(s-t_1)).
\end{eqnarray*}
It then follows by super-additivity that $\lim_{t\to\infty} \frac{1}{t}\ln Z_{\gamma, t}$ exists, although the rate is zero in dimensions $1$ and $2$.

Below we sketch a proof of the precise sub-exponential rates of decay of $Z_{\gamma, t}$ in dimensions $1$ and $2$. We refer ther reader to \cite{DGRS12} for details and proof of the case when $d=3$..

\bigskip

\noindent
{\bf Lower bound on $Z_{\gamma, t}$ for $d=1, 2$:} A lower bound is achieved by clearing a ball of radius $R_t$ around the origin where there are no traps up to time $t$, and then force the random walk $X$ to stay inside this ball up to time $t$. Optimizing the choice of $R_t$ then gives the desired lower bound which, surprisingly, turns out to be sharp in dimensions $1$ and $2$.

Let $B_{R_t}:=\{x\in\Z^d: \Vert x\Vert_\infty\leq R_t\}$ denote the $L_\infty$ ball of radius $R_t$ centered around the origin, with $R_t$ to be optimized over later. Let $E_t$ denote the event that there are no traps in $B_{R_t}$ at time $0$, $F_t$ the event that no traps starting from outside $B_{R_t}$ at time $0$ will enter $B_{R_t}$ before time $t$, and $G_t$ the event that the random walk $X$ with $X(0)=0$ does not leave $B_{R_t}$ before time $t$. Then we have
\begin{equation}\label{315}
Z_{\gamma, t} \geq \P(E_t \cap F_t \cap G_{t}) = \P(E_t)\P(F_t)\P(G_{t}).
\end{equation}

Since $E_t$ is the event that $\xi(0,x)=0$ for all $x\in B_{R_t}$, where $\xi(0, x)$ are i.i.d.\ Poisson with mean $\nu$, we have
\begin{equation}
\P(E_t) = e^{-\nu (2R_t+1)^d}.
\end{equation}
To estimate $\P(G_{t})$, note that for $1\ll R_t\ll \sqrt{t}$, we can approximate $X$ by a Brownian motion and rescale space-time by $(1/R_t, 1/R_t^2)$ to obtain the estimate
\begin{equation}
\P(G_{t}) \geq e^{-c_1t/R_t^2}
\end{equation}
for some $c_1>0$.

To estimate $\P(F_t)$, note that $F_t$ is the event that for each $y\notin B_{R_t}$, no trap starting at $y$ will enter $B_{R_t}$ before time $t$. Since the number of traps starting from each $y\in\Z^d$ are i.i.d.\ Poisson with mean $\nu$, we obtain
\begin{equation}
\P(F_t) = \exp\Big\{-\nu \sum_{y\notin B_{R_t}} \P^Y_y( \tau_{B_{R_t}}\leq t) \Big\},
\end{equation}
where $\tau_{B_{R_t}}$ is the stopping time when $Y$ enters $B_{R_t}$.

In dimension $1$, since we have assumed for simplicity that $Y$ makes nearest-neighbor jumps, we note that $\P(F_t)$ in fact does not depend on the choice of $R_t$. Furthermore, when $R_t=0$, $\P(F_t)$ is easily seen as the annealed survival probability for the trapping problem with instant killing by traps ($\gamma=\infty$) and immobile $X$ ($\kappa=0$). The asymptotics of the annealed survival probability was obtained in \cite[Section 2.2]{DGRS12}, with
\begin{equation}
\P(F_t) = e^{-(1+o(1))\nu \sqrt{8\rho t/\pi}}.
\end{equation}
We can then combine the estimates for $\P(E_t)$, $\P(F_t)$ and $\P(G_t)$ and optimize over $R_t$ to obtain the lower bound
\begin{equation}\label{320}
Z_{\gamma, t} \geq e^{-(1+o(1))\nu \sqrt{8\rho t/\pi}} e^{- ct^{1/3}},
\end{equation}
where the leading order asymptotics is determined by that of $\P(F_t)$, and the second order asymptotics comes from $\P(E_t)\P(G_t)$, with the optimal choice of $R_t$ being a constant multiple of $t^{1/3}$. This lower bound strategy strongly suggests that the fluctuation of the random walk $X$ under the path measure $P_{\gamma, t}$ will be of the order $t^{1/3}$.

In dimension $2$, it was shown in \cite[Section 2.3]{DGRS12} that if $R_t\ll t^\epsilon$ for all $\epsilon>0$, then $\P(F_t)$ has the same
leading order asymptotics as the annealed survival probability of the trapping problem with $\gamma=\infty$ and $\kappa=0$, which coincides with
the asymptotics for $Z_{\gamma, t}$ stated in \eqref{annealed}. To obtain the desired lower bound on $Z_{\gamma, t}$, we can then choose any $R_t$ satisfying $\sqrt{\ln t} \ll R_t\ll t^\epsilon$ for any $\epsilon>0$, which ensures that $\P(E_t)\P(G_t)$ gives lower order contributions. This suggests that under the path measure $P_{\gamma, t}$, $X$ fluctuates on a scale between $\sqrt{\ln t}$ and $t^\epsilon$ for any $\epsilon>0$. However, to identify the optimal choice of $R_t$, we would need to obtain more precise estimates on how $\P(F_t)$ depends on $R_t$.
\bigskip

\noindent
{\bf Upper bound on $Z_{\gamma, t}$ for $d=1, 2$:} The key ingredient is what has been named in the physics literature as the {\em Pascal principle}, which asserts that in \eqref{urep}, if we condition on the random walk trajectory $X$, then the annealed survival probability
\begin{equation}\label{eq:ZXt}
Z^X_{\gamma, t} := \E^\xi\!\left[\exp\!\left\{\!-\gamma\!\int_0^t \!\!\xi(s, X(s))\,{\rm d}s\!\right\}\! \right]\! = \exp \Big\{\nu\!\! \sum_{y\in\Z^d} (v_X(t,y)-1)\Big\}
\end{equation}
is maximized when $X\equiv 0$, provided that the trap motion $Y$ is a symmetric random walk. Therefore
\begin{equation}
Z_{\gamma, t} = \E^X_0[Z^X_{\gamma, t}] \leq Z^{X\equiv 0}_{\gamma, t},
\end{equation}
which has the same leading order asymptotic decay as the case $\gamma=\infty$ and $\kappa=0$ that appeared in the lower bound.

For discrete time random walks under suitable symmetry assumptions, the Pascal principle was proved by Moreau, Oshanin, B\'enichou and Coppey in~\cite{MOBC03, MOBC04}, which can then be easily extended to general continuous time symmetric random walks~\cite{DGRS12}. An interesting corollary is that if $Y$ is a continuous time symmetric random walk, then for any deterministic path $X$ on $\Z^d$ with locally finitely many jumps, we have
\begin{equation}\label{pascal}
\E^Y_0[ |{\rm Range}_{s\in [0,t]}(Y(s))| ] \leq  \E^Y_0[ |{\rm Range}_{s\in [0,t]}(Y(s)+X(s))| ].
\end{equation}
I.e., the expected range of a symmetric random walk can only be increased under deterministic perturbations.

\subsection{Path measures}\label{S:pathmeas}
The only known result so far on the path measures is the following sub-diffusive bound under the annealed path measure $P_{\gamma, t}$ in dimension one, recently proved in \cite{ADS16}.

\begin{theorem}[Sub-diffusivity in dimension one] \label{T:subdiff}
Let $X$ and the trap motion $Y$ follow continuous time random walks on $\Z$ with jump rates $\kappa, \rho>0$ and non-degenerate jump kernels $p_X$ and $p_Y$ respectively. Assume that $p_X$ has mean zero and $p_Y$ is symmetric, with
\begin{equation}
\sum_{x\in\Z} e^{\lambda^* |x|}p_X(x)<\infty \qquad \mbox{and} \qquad \sum_{x\in\Z} e^{\lambda^* |x|}p_Y(x)<\infty
\end{equation}
for some $\lambda^*>0$. Then there exists $\alpha>0$ such that for all $\epsilon>0$,
\begin{equation}\label{pathbdd}
P_{\gamma, t} \big(\Vert X\Vert_t \in (\alpha t^{\frac13}, t^{\frac{11}{24}+\epsilon})\big) \to 1 \quad \text{ as } t \to \infty.
\end{equation}
where $\Vert X\Vert_t:= \sup_{s\in [0,t]} |X_s|$.
\end{theorem}
\noindent
Since $\frac{11}{24}<\frac{1}{2}$, this result shows that under the annealed path measure $P_{\gamma, t}$, $X$ is sub-diffusive.

\begin{remark}
Very recently, \"Oz~\cite{O19} improved the upper bound  $t^{\frac{11}{24}+\epsilon}$ in \eqref{pathbdd} to $ct^{\frac{5}{11}}$ for the continuum analogue model of a Brownian motion among a Poisson field of moving traps.
\end{remark}

We sketch below the proof strategy followed in \cite{ADS16}.
\medskip

\noindent
{\bf Simple random walk, $\gamma=\infty$:} Let us assume for simplicity that $X$ and $Y$ are continuous time simple symmetric random walks, and $\gamma=\infty$. By \eqref{eq:annealedGibbs},
\begin{equation}\label{Pgt1}
P_{\infty, t}(X\in \cdot) = \frac{\E^X_0\big[Z^X_{\infty, t} \, 1_{X\in \cdot}\big]}{Z_{\infty, t}},
\end{equation}
where $Z_{\infty, t} = \E^X_0[Z^X_{\infty, t}]$ with $Z^X_{\gamma, t}$ defined in \eqref{eq:ZXt}, and by \eqref{urep} and \eqref{urep'},
\begin{equation}\label{ZXinftyt}
Z^X_{\infty, t} = \exp \Big\{- \nu \E^Y_0\big[\big|{\rm Range}_{s\in [0,t]}(Y(s)-X(s))\big|\big]\Big\}.
\end{equation}
As shown in the lower bound for the annealed survival probability in \eqref{315}--\eqref{320},
\begin{equation}
Z_{\infty, t} \geq Z^{X\equiv 0}_{\infty, t} e^{-ct^{1/3}}.
\end{equation}
Therefore, from \eqref{Pgt1} we obtain
\begin{eqnarray}
P_{\infty, t}(X\in \cdot) &\leq&  e^{c t^{1/3}} \E^X_0\big[Z^X_{\infty, t}/Z^{X\equiv 0}_{\infty, t} \, 1_{X\in \cdot}\big] \nonumber\\
&=& e^{c t^{1/3}} \E^X_0\Big[e^{- \nu \big(\E^Y_0\big[\big|{\rm Range}_{s\in [0,t]}(Y(s)-X(s))\big| - \big|{\rm Range}_{s\in [0,t]}Y(s)\big|\big]\big)} \, 1_{X\in \cdot}\Big]. \label{Pgt2}
\end{eqnarray}
Note that by the corollary of the Pascal principle \eqref{pascal}, the exponent in the exponential is negative, and hence
\begin{equation}
P_{\infty, t}(X\in \cdot) \leq e^{ct^{1/3}} \P^X_0(X\in \cdot).
\end{equation}
This implies that
\begin{equation}
P_{\infty, t}(\Vert X\Vert_t \leq \alpha t^{1/3}) \leq e^{ct^{1/3}}\P^X_0(\Vert X\Vert_t \leq \alpha t^{1/3}) \to 0 \quad \text{ as } \to \infty,
\end{equation}
if $\alpha$ is sufficiently small, as can be easily seen if $X$ is replaced by a Brownian motion.
This proves the lower bound on the fluctuations.

Proving the sub-diffusive upper bound on $X$ in \eqref{pathbdd} requires finding lower bounds on the differences of the ranges in \eqref{Pgt2} for typical realizations of $X$. Using the fact that $X$ and $Y$ make nearest-neighbor jumps and $Y$ is symmetric, we have
\begin{equation}
\begin{aligned}
\E^Y_0\big[\big|{\rm Range}_{s\in [0,t]}(Y(s)-X(s))\big|\big] & =  1+ \E^Y_0\big[\sup_{s\in [0,t]}(Y(s)-X(s)) - \inf_{s\in [0,t]}(Y(s)-X(s))\big] \\
&= 1+ \E^Y_0\big[\sup_{s\in [0,t]}(Y(s)-X(s))+ \sup_{s\in [0,t]}(X(s)-Y(s))\big] \\
&= 1+ \E^Y_0\big[\sup_{s\in [0,t]}(Y(s)-X(s))+ \sup_{s\in [0,t]}(Y(s)+X(s))\big].
\end{aligned}
\end{equation}
Therefore
\begin{equation}\label{rangedif}
\begin{aligned}
& \E^Y_0\big[\big|{\rm Range}_{s\in [0,t]}(Y(s)-X(s))\big| - \big|{\rm Range}_{s\in [0,t]}Y(s)\big|\big] \\
= \ & \E^Y_0\big[\sup_{s\in [0,t]}(Y(s)-X(s))+ \sup_{s\in [0,t]}(Y(s)+X(s)) - 2 \sup_{s\in [0,t]} Y(s) \big].
\end{aligned}
\end{equation}
We will show that for any $\epsilon>0$, uniformly in $X$ with $\Vert X\Vert_t\geq t^{\frac{11}{24}+\epsilon}$,
\begin{equation} \label{eq:LB}
\text{ the above expectation is bounded from below by $C t^{\frac{1}{3}+\epsilon}$,}
\end{equation}
 which by \eqref{Pgt2} then implies
\begin{equation}\label{334}
P_{\infty, t}(\Vert X\Vert_t \geq t^{\frac{11}{24}+\epsilon})  \leq e^{ct^{\frac13} -Ct^{{\frac13}+\epsilon}} \to 0 \quad \text{ as } t \to \infty.
\end{equation}

To bound the expectation in \eqref{rangedif}, first note that we always have
$$
\sup_{s\in [0,t]}(Y(s)-X(s))+ \sup_{s\in [0,t]}(Y(s)+X(s)) - 2 \sup_{s\in [0,t]} Y(s) \geq 0,
$$
in contrast to $|{\rm Range}_{s\in [0,t]}(Y(s)-X(s))| - |{\rm Range}_{s\in [0,t]}Y(s)|$. Therefore we can restrict to a suitable subset of trajectories of $Y$ to obtain a lower bound.

Let $\sigma_X$ be the first time when $X$ achieves its global maximum in $[0,t]$, and $\tau_X\in [0, t]$ the first time when $X$ achieves its global minimum in $[0,t]$. If $\Vert X\Vert_t\geq t^{\frac{11}{24}+\epsilon}$, then one of the two sets
$$
S:= \{s\in [0, t]: X(\sigma_X) - X(s) \geq t^{\frac{11}{24}+\epsilon}/2\} \quad \mbox{and} \quad T:= \{s\in [0, t]: X(s) - X(\tau_X) \geq t^{\frac{11}{24}+\epsilon}/2\}
$$
must have Lebesgue measure at least $t/2$. Assume w.l.o.g.\ that $|S|\geq t/2$. Then a lower bound on \eqref{rangedif} can be obtained by
restricting $Y$ to the set of trajectories
$$
{\cal F}:= \Big\{\sigma_Y \in S \cap [t/8, 7t/8] \quad \mbox{and} \quad Y(\sigma_Y)-Y(\sigma_X)\leq t^{\frac{11}{24}+\epsilon}/4 \Big\}.
$$
Indeed, on this event,
\begin{equation}\label{335}
\begin{aligned}
& \sup_{s\in [0,t]}(Y(s)-X(s))+ \sup_{s\in [0,t]}(Y(s)+X(s)) - 2 \sup_{s\in [0,t]} Y(s)  \\
\geq \ \ & Y(\sigma_Y)-X(\sigma_Y) + Y(\sigma_X) + X(\sigma_X) - 2 Y(\sigma_Y)  \\
= \ \ & (X(\sigma_X)-X(\sigma_Y)) - (Y(\sigma_Y) - Y(\sigma_X)) \\
\geq \ \ & \frac{1}{4} t^{\frac{11}{24}+\epsilon}.
\end{aligned}
\end{equation}
On the other hand, it is easily seen that there exists $\alpha>0$ such that uniformly in $t$ large,
\begin{equation}\label{336}
\P^Y_0(\sigma_Y \in S \cap [t/8, 7t/8]) \geq \alpha.
\end{equation}
Furthermore, conditioned on $\sigma_Y$, $(Y(\sigma_Y)-Y(\sigma_Y+s))_{s\in [0, t-\sigma_Y]}$ and $(Y(\sigma_Y)-Y(\sigma_Y-s))_{s\in [0, \sigma_Y]}$ are two independent random walks conditioned respectively to not hit zero or go below zero. Such conditioned random walks are comparable to $3$-dimensional Bessel processes. Assume w.l.o.g.\ that $\sigma_X>\sigma_Y$, then we have the rough estimate
\begin{equation}\label{337}
\begin{aligned}
\P^Y_0(Y(\sigma_Y)-Y(\sigma_X) \leq t^{\frac{11}{24}+\epsilon}/4 \, | \sigma_Y) & \geq C \, \P^B_0(\vert B(\sigma_X-\sigma_Y)\vert \leq t^{\frac{11}{24}+\epsilon}/4) \\
& \geq C' \frac{(t^{\frac{11}{24}+\epsilon})^3}{t^{\frac{3}{2}}} = C' t^{-\frac{1}{8}+3\epsilon},
\end{aligned}
\end{equation}
where $B$ is a $3$-dimensional Brownian motion starting from $0$ and its Euclidean norm, $|B_t|$, is a $3$-dimensional Bessel process, while in the last inequality, we used the local central limit theorem for $B$ and $\sigma_X-\sigma_Y\leq t$. Combining \eqref{335}--\eqref{337}, we then obtain
\begin{equation}
\begin{aligned}
& \E^Y_0\big[\sup_{s\in [0,t]}(Y(s)-X(s))+ \sup_{s\in [0,t]}(Y(s)+X(s)) - 2 \sup_{s\in [0,t]} Y(s) \big] \\
\geq\ \ & C' \alpha t^{-\frac{1}{8}+3\epsilon} \frac{t^{\frac{11}{24}+\epsilon}}{4} = C'' t^{\frac{1}{3}+4\epsilon},
\end{aligned}
\end{equation}
which then implies \eqref{eq:LB} and hence \eqref{334}. The above heuristic calculations were made rigorous in \cite{ADS16}.
\bigskip

\noindent
{\bf Simple random walk, $\gamma<\infty$:} When $\gamma<\infty$, the representation \eqref{ZXinftyt} for $Z^X_{\infty, t}$ should be replaced by
\begin{equation}\label{ZXgt}
Z^X_{\gamma, t} =  \exp \Big\{- \nu \E^{Y, \cal T}_0\big[\big|{\rm SoftRange}_{s\in [0,t]}(Y(s)-X(s))\big|\big]\Big\}
\end{equation}
as can be seen from \eqref{Zrep2}. The difficulty then lies in controlling the difference
\begin{equation}
\begin{aligned}
F_{\gamma, t}(Y-X) & :=  \big[\big|{\rm Range}_{s\in [0,t]}(Y(s)-X(s))\big|\big] - \E^{\cal T}\big[\big|{\rm SoftRange}_{s\in [0,t]}(Y(s)-X(s))\big|\big]\\
& = \sum_{x\in \Z} e^{-\gamma L_t^{Y-X}(x)} 1_{L_t^{Y-X}(x)>0},
\end{aligned}
\end{equation}
where $L_t^{Y-X}(x) := \int_0^t 1_{\{Y(s)-X(s)=x\}} {\rm d}s$ is the local time of $Y-X$ at $x$. In \cite{ADS16}, this control is achieved by proving that
\begin{equation}\label{341}
\sup_{t\geq e} \E^{Y}_0\Big[\exp\Big\{\frac{c}{\ln t} F_{\gamma, t}(Y)\Big\}\Big] <\infty,
\end{equation}
which also leads to interesting bounds on the set of thin points of $Y$, i.e., the set of $x$ where the local time $L^Y_t(x)$ is positive but unusually small.
\bigskip

\noindent
{\bf Non-simple random walks:} When $X$ and $Y$ satisfy the assumptions in Theorem~\ref{T:subdiff}, but are not necessarily simple random walks,
the arguments outlined above can still be salvaged, provided we can control the difference
\begin{equation}
\begin{aligned}
G_t(Y) & := \sup_{s\in [0,t]} Y(s) - \inf_{s\in [0,t]}Y(s) + 1 - |{\rm Range}_{s\in [0,t]} (Y(s))| = \sum_{\underset{s \in [0,t]}{\inf} Y(s)\leq x\leq \underset{s \in [0,t]}{\sup} Y(s)} 1_{L^Y_t(x)=0},
\end{aligned}
\end{equation}
which is the total size of the holes in the range of $Y$. In \cite{ADS16}, this control is achieved by proving
\begin{equation}\label{343}
\sup_{t\geq e} \E^{Y}_0\Big[\exp\Big\{\frac{c}{\ln t} G_t(Y)\Big\}\Big] <\infty.
\end{equation}
The proof of \eqref{341} and \eqref{343} in \cite{ADS16} in fact follow the same line of arguments.

\section{Some Open Questions}\label{S:open}

The large deviation and spectral techniques that have been successful for the trapping problem with immobile traps largely fail for mobile traps. As a result, many questions remain unanswered for the trapping problem with mobile traps, even when the traps are just a Poisson system of independent random walks. We list below some natural open questions.
\medskip

\noindent
{\bf Open questions:}

\begin{enumerate} [label=(\arabic*)]
\item In dimension $1$, we conjecture that under the annealed path measure $P_{\gamma, t}$, $X$ fluctuates on the scale $t^{1/3}$, which is based on the lower bound strategy for survival in Section \ref{S:survprob}. In fact, we saw that the probability that traps from outside the ball of radius $R_t$ do not enter it before time $t$ does not depend on the choice of $R_t$, which leaves only the interplay between the cost of clearing the ball of traps at time $0$ and the cost of the forcing the walk to stay within the ball up to time $t$. This is also what happens in the case of immobile traps, which leads us to conjecture that not only the fluctuation is on the same scale of $t^{1/3}$~\cite{S90, S03}, but also the rescaled paths converge to the same limit.

\item In dimension $2$, the lower bound strategy for the annealed survival probability suggests that under the annealed path measure $P_{\gamma, t}$, $X$ fluctuates on a scale $R_t$ with $\sqrt{\ln t} \ll R_t \ll t^{\epsilon}$ for all $\epsilon>0$. The correct scale of $R_t$ still needs to be determined. In dimension $1$, we heavily used the fact that the range of a random walk is essentially captured by its maximum and minimum. This is no longer applicable in dimensions $d\geq 2$, and new techniques need to be developed.

\item In dimensions $d\geq 3$, if we follow the same lower bound strategy for the annealed survival probability, then it is not difficult to see that the optimal size of the ball which is free of traps should be of order one instead of diverging as in dimensions $1$ and $2$. This seems to suggest that under the annealed path measure $P_{\gamma, t}$, $X$ should be localized near the origin. On the other hand, it is also reasonable to expect that the interaction between the traps and the random walk $X$ will create a clearing around $X$ in a time of order $1$, and random fluctuations will then cause the clearing and $X$ to undergo diffusive motion, leading to a central limit theorem and invariance principle for $X$ under $P_{\gamma,t}$. We conjecture that the second scenario is what actually happens.

\item What can we say about the quenched path measure $P^\xi_{\gamma, t}$? For Brownian motion among Poisson obstacles, Sznitman~\cite{S98} has shown that the particle $X$ seeks out pockets of space free of traps. There is a balance between the cost of going further away from the starting point to find larger pockets and the benefit of surviving in a larger pocket, which leads to super-diffusive motion under the quenched path measure (as pointed out to us by R. Fukushima, this follows from results in \cite{S98}; and see \cite{DX17} for the Bernoulli hard trap model). When the traps are mobile, these pockets of space free of traps are destroyed quickly and become much more rare. Would $X$ still be super-diffusive under the quenched path measure?

\item In our analysis, we heavily used the fact that the random trapping potential $\xi$ is a Poisson field. One could also consider $\xi$ generated by other interacting particle systems, such as the exclusion process or the voter model as G\"artner et al considered in the context of the parabolic Anderson model~\cite{GHM09, GHM12}. Their analysis has focused on the exponential asymptotics of the quenched and annealed solution of the PAM. The path behavior of the corresponding trapping problem remains open. The behavior is expected to be similar when $\xi$ is the occupation field of either the symmetric exclusion process or the Poisson system of independent random walks, because the two particle systems have the same large scale space-time fluctuations.

\item It is natural to also consider the case when the random walk
  moves with a deterministic drift. A naive survival
  strategy of creating a ball centered at the origin free of traps,
  then there is an exponential cost for the random walk to stay within
  that region, while leaving the region will also incur an exponential
  cost depending on $\gamma$, the rate of killing when the walk meets
  a trap. Depending the strength of the drift relative to $\gamma$, it
  is conceivable that there is a transition in the path behavior as
  the drift varies, such that when the drift is small, the random walk
  is localized near the origin just as in the case of zero drift, and
  when the drift becomes sufficiently large, the random walk becomes
  delocalized. This is also a question of active interest for the PAM
  with immobile traps (see e.g.~\cite[Sec.~7.10]{K16}, \cite[Sec.~5.4
    \& 7.3]{S98}).

\item In Theorem \ref{T:annealed} on the asymptotics of the annealed survival probability, we note that in dimensions $1$ and $2$, the precise rate of decay of $Z_{\gamma, t}$ does not depend on the killing rate $\gamma$. The heuristics is that as long as $\gamma>0$, the random walk $X$ would not see the traps under the annealed path measure. This suggests sending $\gamma=\gamma_t\to 0$ as $t\to\infty$ in order to see non-trivial dependence on the killing rate. In dimension $1$, the correct rate is $\gamma_t=\widetilde \gamma/\sqrt{t}$, and in dimension $2$, $\gamma_t=\widetilde \gamma/\ln t$. Decaying killing rate has been considered in the immobile trap case, see \cite[Sec.~7.3.3]{K16} and the references therein.
\end{enumerate}

\bigskip

\noindent
{\bf Acknowledgement.} R.S. is supported by NUS grant R-146-000-220-112. S.A. is supported by CPDA grant and ISF-UGC project. We thank Ryoki Fukishima for helpful comments that corrected some earlier misstatements. Lastly, we thank Vladas Sidoravicius for encouraging us to write this review, and we are deeply saddened by his untimely death.

\end{document}